\documentclass{ws-ijbc}
\usepackage{ws-rotating}     
\usepackage{graphicx}
\usepackage{epstopdf}
\usepackage{amsmath}
\usepackage{amssymb}

\newcommand{\D}{\mbox{\bf D}}

\newcommand{\Z}{\mbox{$\mathbb{Z}$}}
\newcommand{\AND}{\quad\mbox{and}\quad}
\newcommand{\Otwo}{{\bf O(2)}}
\newcommand{\R}{\mbox{$\mathbb{R}$}}
\newcommand{\N}{\mbox{$\mathbb{N}$}}
\newcommand{\Fix}{\mbox{{\rm Fix}}}
\newcommand{\qed}{\hfill\mbox{\raggedright\rule{0.07in}{0.1in}}\vspace{0.1in}}

\begin{document}

\catchline{}{}{}{}{} 

\markboth{Buono {\sl et al.}}{Hidden symmetry in a Kuramoto-Sivashinsky IBVP}

\title{Hidden symmetry in a Kuramoto-Sivashinsky initial-boundary value problem}

\author{Pietro-Luciano Buono, Lennaert van Veen and Eryn Frawley\footnote{Corresponding author: Eryn.Frawley@uoit.net}}

\address{Faculty of Science, University of Ontario Institute of Technology, 2000 Simcoe Street North\\
Oshawa, ON L1H7K4, Canada}

\maketitle

\begin{history}
{\it Under review}
\end{history}

\begin{abstract}
We investigate the bifurcation structure of the Kuramoto-Sivashinsky equation with homogeneous Dirichlet boundary conditions. Using hidden symmetry principles, based on
an extended problem with periodic boundary conditions and $\Otwo$ symmetry, we show that the zero solution exhibits two kinds of pitchfork bifurcations: one that breaks
the reflection symmetry of the system with Dirichlet boundary conditions and one that breaks a shift-reflect symmetry of the extended system. Using Lyapunov-Schmidt reduction,
we show both to be supercritical. We extend the primary branches by means of numerical continuation, and show that they lose stability in pitchfork, transcritical or Hopf bifurcations.
Tracking the corresponding secondary branches reveals an interval of the viscosity parameter in which up to four stable equilibria and time-periodic solutions coexist.
Since the study of the extended problem is indispensible for the explanation of the bifurcation structure, the Kuramoto-Sivashinsky problem with Dirichlet boundary conditions
provides an elegant manifestation of hidden symmetry.
\end{abstract}


\section{Introduction}

In this paper, we explore the bifurcation structure of the Kuramoto-Sivashinsky (KS) equation in one space dimension with Dirichlet Boundary Conditions (BC) as a function of the viscosity parameter. In order to give a fairly complete overview of the bifurcation behaviour close to the leading instability, we use three techniques: an equivariant branching lemma, Lyapunov-Schmidt reduction and numerical continuation of equilibria and time-periodic orbits.

The KS equation was derived from work by \citet{Kuramoto} on the Belousov-Zhabotinskii equation and \citet{Siva} on the dynamics of laminar flame fronts. For a summary of the origin, derivation and importance of the KS equation see, e.g.,~\citet{Misbah-Valence-1994} or \citet{WH-1999}. Here, we pose it in the following form:
\begin{align}\label{KS}
u_t+u u_x + u_{xx}+\nu u_{xxxx}&=0.
\end{align}
The Dirichlet boundary conditions are
\begin{align}\label{KS-DBC}
u(-1,t)=u_{xx}(-1,t)=u(1,t)=u_{xx}(1,t)&=0.
\end{align}
Equation (\ref{KS}) has a rich bifurcation structure, which has been investigated by several authors. Most studies have focused on the case of periodic BC, e.g.~\citet{WH-1999,Cvitanovic-etal-2010,Zhang-etal-2011}. In this setting, the KS equations is $\Otwo$ symmetric and the consequences of this symmetry have been explored in detail ~\cite{Cvitanovic-etal-2010}. There are fewer studies of the KS equation with Dirichlet BC (\ref{KS-DBC}). The existence and uniqueness of solutions was established by \citet{Galak}. \citet{LiChen} studied bifurcations from the zero solution with Dirichlet, Neumann and periodic BC. They identify families of bifurcation points for each case and use Lyapunov-Schmidt reduction to determine the criticality of the bifurcating solutions. For the Dirichlet BC case, they identify a family of steady-state bifurcation points with a kernel generated by $\sin(n \pi x)$ for $n\in \N$.

Our first result is the existence of a second family of bifurcation points, with kernel $\cos ([2n-1]\pi x/2)$ for $n\in \N$.
This includes the leading bifurcation point, i.e. the first instability to occur as the viscosity parameter is reduced from unity, where the zero solution is stable.
The KS initial-boundary value problem with Dirichlet BC has a $\mathbb{Z}_2$ symmetry, and the bifurcations with kernel $\cos([2n-1] \pi x/2)$ can readily be explained as
breaking this symmetry. Those with kernel $\sin(n\pi x)$ give rise to branches of solutions that lie entirely in the fixed point subspace of the symmetry. Their presence can be explained
if we exploit the fact that solutions of the Dirichlet BC case can be seen as restrictions of solutions of the periodic BC case on an extended domain;
this is the phenomenon of hidden symmetry. The extended system has $\Otwo$ symmetry and, accordingly, two-dimensional kernels. Using the structure of the extended problem, we formulate
an equivariant branching lemma that determines all bifurcations of the zero state to be of the pitchfork type.
In addition, we use Lyapunov-Schmidt reduction to show that all these pitchfork bifurcations are supercritical.

In order to study the nonlinear behaviour of the system, we extend the primary branches by means of numerical continuation. These, too, exhibit
pitchfork bifurcations, in addition to transcritical bifurcations, that are neatly explained in terms of the hidden symmetry. All bifurcation points at which primary branches
lose stability are identified and the corresponding secondary branches are computed, including two branches of time-periodic solutions. Thus, we are able to give a fairly complete
overview of the dynamics of the KS initial-boundary value problem around the first four instabilities of the zero solution.

The paper is organized as follows. Section~\ref{sec:lin} discusses the linearization at the zero solution and we determine the bifurcation values for the Dirichlet and periodic BC cases. In Section~\ref{sec:hidden}, we embed the KS equation with Dirichlet BC on $(-1,1)$ in a periodic BC problem on the doubled-up interval from $(-2,2)$ and derive results from hidden symmetry theory~\cite{GolSte}.
This is followed in Section~\ref{sec:LSred} by the computation of the Lyapunov-Schmidt reduction for the bifurcation points
at which the kernel is given by $\cos([2n-1]\pi x/2)$. Finally, Section~\ref{numerical} shows the result of numerical continuations of primary and secondary
branches of equilibria and time-periodic solutions.

\section{Linearization and Kernels}\label{sec:lin}
It is straightforward to verify that the function $u^{*}\equiv 0$ on $(-1,1)$ is a homogeneous steady-state solution of~(\ref{KS}) satisfying~(\ref{KS-DBC}).
The linearization of~(\ref{KS}) at $u^{*}$ is
\[
{\cal L}(v) := v_t + v_{xx}+\nu v_{xxxx}=0.
\]
At the homogeneous equilibrium $u^{*}=0$, writing $v(x,t) = (a\cos(k\pi x)+b\sin(k\pi x))e^{\lambda t}$,
the characteristic equation is
\[
\lambda - k^2\pi^2 + \nu k^4\pi^4 = 0
\]
and for $\lambda=0$ we have $\nu = (k\pi)^{-2}$. The boundary conditions for $v$ are satisfied for $a=0,\,k=n$ for $n\in \N$ as shown
by~\citet{LiChen}, but also for $b=0,\,k = n-\frac{1}{2}$ where $n\in \N$. The smallest positive value of $k$ is $1/2$
and therefore the largest value of $\nu$ at which a bifurcation point arises is
\[
\nu^{*} = \frac{4}{\pi^2}.
\]
The kernel of ${\cal L}$ is one-dimensional at each bifurcation point. For $k=n$, the kernel is spanned by $e^{\rm o}_n(x)=\sin(n\pi x)$
while for $k=n-\frac{1}{2}$, the kernel is spanned by $e^{\rm e}_n(x)=\cos([2n-1]\pi x/2)$.

We embed the original partial differential equation with Dirichlet boundary conditions on $(-1,1)$
into the domain $(-2,2)$ and we impose periodic boundary conditions
\begin{equation}\label{eq:per-BC}
\dfrac{\partial^{\ell} }{\partial x^{\ell}}u(-2,t)=\dfrac{\partial^{\ell} }{\partial x^{\ell}}u(2,t), \qquad \ell=0,1,2,3.
\end{equation}
For the extended system, at the homogeneous equilibrium $u^{*}=0$, the characteristic equation is unchanged from above
and for $\lambda=0$ we still have $\nu = (k\pi)^{-2}$. We obtain the same bifurcation values as the Dirichlet BC case on
the interval $(-1,1)$. The eigenspace is two-dimensional and one can directly check that the kernel is generated by the functions
$\{\cos k\pi x, \sin k\pi x\}$. We summarize the results in the following statement.

\begin{proposition}\label{prop:bif}
For the Dirichlet and periodic BC cases, the linearization ${\cal L}$ at the homogeneous zero solution has bifurcation
points at $\nu^{*}=(k\pi)^{-2}$ with $k=n$ and $k=n-\frac{1}{2}$, $n\in \N$. In the Dirichlet BC case, for $k=n$, the kernel is spanned by $e^{\rm o}_n(x)$
while for $k=n-\frac{1}{2}$, the kernel is spanned by $e^{\rm e}_n(x)$. In the periodic BC case, the kernel is spanned by
$\{\cos \pi k x,\sin \pi k x\}$ both for $k=n$ and $k=n-\frac{1}{2}$.
\end{proposition}

In the next section, we explore the relationship between the solutions of the Dirichlet BC problem on $(-1,1)$ and the periodic BC case on $(-2,2)$
defined above.

\section{Hidden symmetry}\label{sec:hidden}
We now discuss the hidden symmetry properties of~(\ref{KS}) with Dirichlet boundary conditions~(\ref{KS-DBC}) and in particular, show
that all steady-state bifurcations from the trivial solution must be pitchfork bifurcations.

It is well-known~\cite{GolSte}, that an equation such as~(\ref{KS}) with boundary conditions~(\ref{eq:per-BC}) is symmetric
with respect to a group generated by reflections and translations. In particular, if $u(x,t)$ is a solution then
\begin{equation}\label{eq1}
\kappa.u(x,t) = -u(-x,t) \AND \theta.u(x,t) = u(x-\theta,t), \quad \theta\in (-2,2]
\end{equation}
are also solutions. In particular, $\theta^{-1}\kappa=\kappa\theta$ for all $\theta\in (-2,2]$.
The group generated by $\langle \kappa,\theta\rangle$ is isomorphic to the group $\Otwo$.

We drop the dependence on $t$ and define the space
\[
X_{per} = \{ u \in C^4([-2,2],\R)\mid \mbox{$u$ satisfies boundary conditions~(\ref{eq:per-BC})} \}.
\]
The group $\Gamma = \Otwo$ acts on this space of functions by the formulae given in~(\ref{eq1}).

At the bifurcation points $(u^{*},\nu_k^{*}) = (0,(k\pi)^{-2})$ for $k=n-\frac{1}{2}$ and $k=n$ with $n\in \N$, the
kernel of the linearization at $\nu_k^{*}$ is two-dimensional and generated by the functions $\{\cos k\pi x, \sin k\pi x\}$ as noted in
Proposition~\ref{prop:bif}.

The kernel of the linearization is an absolutely irreducible representation of $\Otwo$. From the Equivariant Branching
Lemma~\cite{GSS88}, a branch of equilibrium solutions for each isotropy subgroup with fixed point subspace
of dimension one bifurcates from $(u^{*},\nu_k^{*})$. The two-dimensional $\Otwo$ irreducible representation spanned
by $\{\cos k\pi x, \sin k\pi x\}$ has a family of isotropy subgroups $\Sigma_{k,\theta} = \langle \theta^{-1}\kappa\theta, 2/k\rangle$
with $2/k$ a translation, each isomorphic to $\D_{2k}$, where $\theta \in [-2/k,2/k)$. $\Sigma_{k,\theta}$ fixes the one-dimensional subspace
spanned by $\sin k\pi\theta\cos k\pi x+\cos k\pi\theta\sin k\pi x=\sin(k\pi(x+\theta))$. Note that the translation by $1/k$ acts by $-1$ on
$\Fix(\Sigma_{k,\theta})$ for all $\theta$:
\[
(1/k).\sin(k\pi(x+\theta)) = \sin(k\pi(x+\theta))\cos(\pi) = -  \sin(k\pi(x+\theta)).
\]
Finally, one can verify that $\Fix(\Sigma_{k,1/2k})=\Fix(\Sigma_{k,-1/2k}) = \mbox{span}\{\cos k\pi x\}$ and
\[
\Fix(\Sigma_{k,1/k})=\Fix(\Sigma_{k,-1/k})=\Fix(\Sigma_{k,2/k})=\Fix(\Sigma_{k,0}) = \mbox{span}\{\sin k\pi x\}
\]

\begin{proposition}
The periodic BC problem has $\Otwo$ steady-state bifurcation points at $(0,\nu_{k}^{*})$ where
$\nu_{k}^{*}=(k\pi)^{-2}$ for $k=n-1/2$ and $k=n$ with $n\in \N$. An $\Otwo$-orbit of steady-state solutions
with isotropy subgroup $\Sigma_{k,\theta}$ conjugate to $\D_{2k}$ bifurcates from $(0,\nu_{k}^{*})$ where solutions
with the same isotropy subgroup come in pairs related by the translation element $\theta_{1/k}$.
\end{proposition}

\begin{remark}\label{rem:dihedral}
At bifurcation points $\nu_{k}^{*}$ for $k=n-1/2$, the symmetry group $\Sigma_{k,\theta}$ has cyclic symmetry of order $2n-1$
while for $k=n$ the cyclic symmetry is of order $2n$.
\end{remark}

%

We now present the connections between the solutions of the Dirichlet BC case with the periodic BC case,
We begin with this result concerning the regularity of solutions related by hidden symmetry. The proof describes
how equilibrium solutions are extended from the smaller to the larger domain. The proof follows the one given in~\citet{GolSte}
for reaction-diffusion equations.
\begin{lemma}
Smooth equilibrium solutions of the Dirichlet BC problem on $[-1,1]$ extend to smooth equilibrium solutions
of the periodic BC problem on $[-2,2]$.
\end{lemma}

\proof Let $u(x)$ be an equilibrium solution of~(\ref{KS}) with Dirichlet BC on $[-1,1]$; the time-dependence of
$u$ is suppressed as it is not needed. Extend $u(x)$ to $\hat{u}(x)$ defined by
\[
\hat{u}(x) = \left\{ \begin{array}{cc} -u(-2-x), & x\in [-2,-1]\\
u(x), & x\in [-1,1] \\
-u(2-x), & x\in [1,2].
\end{array} \right.
\]
The function $\hat{u}$ is constructed from $u$ by extending in an odd way with respect to $x=1$
the portion defined over $[0,1]$ to $[1,2]$. Similarly, extending in an odd fashion the portion over
$[-1,0]$ with respect to $x=-1$ to $[-2,-1]$.

Then, $\hat{u}(x)$ satisfies~(\ref{KS}-\ref{KS-DBC}) automatically on $[-1,1]$ and we show the proof for the portion
of $\hat{u}$ defined on $[-2,-1]$. For $x\in [-2,-1]$, let $y=-2-x$, then $y\in [0,1]$ and the partial derivatives
with respect to $x$ become partial derivatives with respect to $y$ so that
\[
\hat{u}(x)\hat{u}_{x}(x)+\hat{u}_{xx}(x)+\nu \hat{u}_{xxxx}(x) = u(y,t)u_{y}(y) + u_{yy}(y) + \nu u_{yyyy}(y) = 0.
\]
An identical computation holds for the interval $[1,2]$. It is straightforward to check that the
periodic boundary conditions are satisfied. The smoothness needs to be verified only at $x=\pm 1$.
We do the case at $x=1$, the computations at $x=-1$ are similar and omitted.
We use the notation $1^{-}$ and $1^{+}$ to denote the left and right limits at $x=1$.
By construction of $\hat{u}$, we have $\hat{u}(1^{-})=\hat{u}(1^{+})$. Note that for $x\in (1,2)$
the $k^{th}$ derivative satisfies $\hat{u}^{(k)}(x) = (-1)^{k+1}u^{(k)}(2-x)$. Thus,
$\hat{u}^{(2k+1)}(1^{+})=u^{(2k+1)}(1^{-})=\hat{u}^{(2k+1)}(1^{-})$ for all $k\in \N$. Now,
$\hat{u}''(1^{+})=-u''(1^{-})=\hat{u}''(1^{-})$, but we know that $u''(1)=0$ because $u$
satisfies the Dirichlet BC. Finally, using~(\ref{KS-DBC}) we have
\[
u^{(4)}(\pm 1) = -\nu^{-1}(u(\pm 1)u_{x}(\pm 1) + u_{xx}(\pm 1)) = 0
\]
and this implies $\hat{u}^{(4)}(1^{-})=\hat{u}^{(4)}(1^{+})=0$. We show by induction that
$u^{(2k)}(\pm 1)=0$ with $k\geq 3$. For $k=3$,
\[
u^{(6)}(\pm 1) = - \nu^{-1}(3u_{x}(\pm 1)u_{xx}(\pm 1) + u_{xxxx}(\pm 1)) = 0.
\]
For $u^{(2k)}(x)$, one can verify that the expression for this derivative depends only
on terms of the form $u^{\ell}u^{n}$ where $\ell$ is odd and $n$ is even and the term
$u_{xx}$ obtained from substituting for $u^{(4)}(x)$. But the second derivative of terms
such as $u^{\ell}u^{n}$ is in terms of a sum of terms of the same form. Thus, evaluating
at $\pm 1$, those all vanish, and we know $u_{xx}(\pm 1)=0$, thus $u^{(2k)}(\pm 1)=0$.
Therefore, $\hat{u}$ is smooth. \qed

We can characterize using the symmetry reflection the solutions of the periodic BC problem which satisfy the Dirichlet BC.
\begin{lemma}\label{lemma:restriction}
Let $v$ be a steady-state solution of the extended periodic BC problem on $[-2,2]$.
Then, $v=\hat{u}$ where $u$ satisfies the Dirichlet BC problem
if and only if $v\in \Fix((1)^{-1}\Z_2(\kappa)(1))$.
\end{lemma}

\proof Suppose $v=\hat{u}$ satisfies the periodic BC problem. Then, by construction of $\hat{u}$,
we have
\[
v(x-1) = \left\{ \begin{array}{cc}
-u(1-x), & x\in [-1,0]\\
u(x-1), & x\in [0,2] \\
-u(1-x), & x\in [-2,-1].
\end{array} \right.
\]
and so $v(x-1)$ is an odd function, that is, $v(x-1)\in \Fix(\Z_2(\kappa))$. This
means $v(x)\in \Fix(1^{-1}\Z_2(\kappa) 1)$. The opposite implication
is straightforward because elements in the subspace $\Fix(1^{-1}\Z_2(\kappa) 1)$
automatically satisfy the Dirichlet BC.  \qed

The main result of this section concerns the steady-state bifurcation points of~(\ref{KS}) and~(\ref{KS-DBC}).
\begin{theorem}\label{thm:DBC-sol}
The Dirichlet BC problem~(\ref{KS}) and (\ref{KS-DBC}) has, generically, pitchfork steady-state bifurcation points at $(0,\nu_{k}^{*})$
for $k=n-1/2$ and $k=n$ for $n\in \N$. If $k=n-1/2$, the pair of bifurcating solutions are on the same $\kappa$ orbit.
If $k=n$, the pair of bifurcating solutions are fixed by $\Z_2(\kappa)$, and therefore not related by symmetry.
\end{theorem}

\proof At the bifurcation points $(0,\nu_k^{*})$ for $k=n-1/2$ in the periodic BC case, there are two branches of steady-state
solutions (a $1/2k$ group orbit) with isotropy subgroup $\Sigma_{k,1/2k}\supset (1)^{-1}\Z_2(\kappa) (1)$.
By Lemma~\ref{lemma:restriction}, the pair of bifurcating solutions also satisfy the Dirichlet BC problem~(\ref{KS})
and~(\ref{KS-DBC}). We now show that for the periodic BC problem, $1/k$ acts on $\Fix(\Sigma_{k,1/2k})$
in the same way as $\kappa$. Let $u(x)\in \Fix(\Sigma_{k,1/2k})$ and recall that $(1/2k).u(x)$ is
fixed by $\kappa$. Then, using the non-commutativity relationship of $\theta$ and $\kappa$,
\[
(\kappa (1/k)).u(x) = (\kappa (1/2k)).u(x-1/2k) =(1/2k)^{-1}\kappa.u(x-1/2k) = (1/2k)^{-1}u(x-1/2k)=u(x)
\]
which indeed means $\theta_{1/k}.u(x)=-u(-x)=\kappa.u(x)$. Therefore, the pair of bifurcating solutions
is related by the $\kappa$ symmetry which is the only remaining symmetry of the Dirichlet BC problem.

Consider the $k=n$ bifurcation points. The kernel at those bifurcating points of the periodic BC problem
restricted to $\Fix(\Sigma_{k,0})$ has a pair of bifurcating branches of steady-state solutions satisfying the
the Dirichlet BC problem. But $\Z_2(\kappa)\subset \Sigma_{k,0}$ and so the bifurcating pair for the Dirichlet
BC problem are not related by a symmetry. \qed.

\begin{remark}
For the cases $k=n\in \Z$ of the Dirichlet BC problem, Li and Chen~\cite{LiChen} show that the cubic coefficient
obtained via Lyapunov-Schmidt reduction is negative therefore proving that the pitchfork bifurcations are
supercritical. It is also shown in~\cite{LiChen} that the periodic BC case for the cases $k=n$ restricted to the fixed point
subspaces with kernels consisting of $\sin(k\pi x)$ and $\cos(k\pi x)$ are also supercritical. As this case has in
fact $\Otwo$-symmetry, it is only necessary to verify one of them to determine the supercritical nature of the
bifurcation of the $\Otwo$ group orbit.
\end{remark}
In the next section, we complete the analysis by using the Lyapunov-Schmidt reduction to determine the
criticality at the bifurcation points with $k=n-1/2$ for the Dirichlet BC case.

We conclude this section by looking at secondary bifurcations, that is, bifurcations from the primary branches
of bifurcation emanating from the $u^{*}=0$ steady-state. We use again the periodic BC case to obtain the results and
we separate the $k=n-1/2$ and $k=n$ cases.

We begin with the $k=n-1/2$ case and recall from Remark~\ref{rem:dihedral} that the branches of steady-states
bifurcating at $\nu_{k}^{*}$ have isotropy subgroup $\Sigma_{k,\theta}\simeq \D_{2n-1}$
Standard results from
equivariant bifurcation theory with dihedral symmetry \cite{GSS88} state that generic symmetry-breaking bifurcations from a steady-state with isotropy subgroup $\D_{2n-1}$ lead to a unique (up to conjugacy) branch of steady-state solutions with isotropy subgroup given by the reflection symmetry. We focus on the isotropy subgroup $\Sigma_{k,1/2k}$ and so the isotropy subgroup of a bifurcating branch from the $\Sigma_{k,1/(2k)}$ branch is $\Z_2((1/(2k))^{-1}\kappa (1/(2k)))$ and there are no other element of $\Otwo$ leaving $\Fix(\Z_2((1/(2k))^{-1}\kappa (1/(2k)))$ invariant; that is, $\Z_2((1/(2k))^{-1}\kappa (1/(2k)))$ is its own normalizer subgroup. From Theorem~\ref{thm:DBC-sol}, to the branch of steady-state solutions with isotropy subgroup $\Sigma_{k,1/(2k)}$ there corresponds a branch of steady-state solutions of the Dirichlet BC problem. Therefore, if the steady-state solution branch with $\Sigma_{k,1/2k}$ isotropy subgroup has a symmetry-breaking steady-state bifurcation to a branch of solutions which also satisfy the Dirichlet BC, then, generically, it appears as a transcritical bifurcation in the Dirichlet BC problem.

The case $k=n$ leads to isotropy subgroups $\Sigma_{k,\theta}\simeq \D_{2n}$ and generic steady-state bifurcation results with dihedral symmetry $\D_{2n}$ leads to pitchfork bifurcations because the normalizer subgroup of the $\Z_2$ symmetric branches acts by $-1$ on the fixed point subspaces. Focusing on the isotropy subgroup $\Sigma_{k,0}$, we know this branch of steady-state solutions also satisfies the Dirichlet BC case and if it has a symmetry-breaking bifurcation which also satisfies the Dirichlet BC, then it appears as a pitchfork bifurcation in the Dirichlet BC problem.


\section{Lyapunov-Schmidt reduction}~\label{sec:LSred}
We now determine the criticality of the primary branches of steady-states bifurcating from the bifurcation points
at $\nu_{k}^{*}$ for $k=n-1/2$. We use the Lyapunov-Schmidt reduction to perform this analysis. For the Lyapunov-Schmidt
reduction process below, it is convenient to rewrite the time-independent part of the equation abstractly as follows. We define
the function space
\[
X=\{u \in C^4[-1,1] \mid u(-1)=u(1)=u''(-1)=u''(1)=0\} ,
\]
and let $Y=C^0[-1,1]$. The two spaces share the inner product defined as
\[
\langle u,v \rangle = \int_{-1}^{1}u(\xi)v(\xi)d\zeta.
\]
\begin{theorem}\label{thm:critical}
For all the bifurcating points $\nu_{k}^{*}$ with $k=n-1/2$ with $n\in \N$, the branches of steady-states
bifurcate supercritically. For the largest bifurcation point $\nu_{1/2}^{*}$, the steady-state solutions bifurcating
are asymptotically stable.
\end{theorem}

\proof We rewrite the steady-state Kuramoto-Sivashinsky equation using the mapping $\Phi:X\times\R\to Y$ as
\begin{align}\label{KS-op}
\Phi(u, \nu) = uu_x + u_{xx}+\nu u_{xxxx}=0.
\end{align}
We define the linear operator
$L_0: X \times \R \mapsto Y$ as follows
\[
L_0 = d\Phi(0,\nu_{k}^{*})=\nu_{k}^{*} \frac{d^4}{d\xi^4}+\frac{d^2}{d\xi^2}
\]
and we remove on purpose the dependence of this linear operator on $k$ to lighten the notation. 
$L_0: X \rightarrow Y$ is self-adjoint because of the boundary conditions defining $X$. Furthermore, $\ker L_0=\mbox{span}\{e\}$, where
$e =e^{\rm e}_n(\xi)$.  To keep the notation as simple as possible below, we let
$\beta = \frac{\pi}{2}(2n-1)$. The subspaces $\ker L_0$ and $\mbox{range}\,L_0$ enable us to decompose the spaces
$X$ and $Y$:
\[
X=\ker L_0 \oplus (\ker L_0)^\bot  \AND   Y =(\mbox{range}\,L_0)^\bot \oplus \mbox{range}\,L_0.
\]
Notice that because $L_0$ is self-adjoint, then $\ker L_0=(\mbox{range}\, L_0)^{\bot}$.
We can now define $P$ as the orthogonal projection from $Y$ onto range $L_0$. Let $u \in X$, then the decomposition of $X$ leads to $u=ye+w$, where $y \in \R$ is the coordinate on the kernel and $w \in (\ker L_0)^{\bot}$. The projection $P$ can be applied to $\Phi$ to decompose it in two separate parts, namely
\[
G(y,w,\nu):= P\Phi(ye+w,\nu)=0 \AND (I-P)\Phi(ye+w,\nu)=0.
\]
Because $G(0,0,\nu_{k}^{*})=0$ and $dG(0,0,\nu_{k}^{*})$ is surjective on $\mbox{range}\,L_0$, by the implicit function theorem, there exists a unique function $w=w(y,\nu)$ with $w(0,\nu_{k}^{*})=0$ such that $G(y,w(y,\nu),\nu)\equiv 0$.
The equation $(I-P)\Phi(ye+w(y,\nu),\nu)=0$ is transformed by taking the inner product with $e$, from which we obtain
\[
g(y,\nu) = \langle e, \Phi(ye + w(y,\nu),\nu)\rangle = 0
\]
since $P\Phi$ is orthogonal to $e$.
We want to find an approximation of the system defined by $g(y,\nu)$ at the point $(0,\nu_{k}^{*})$, so it is necessary to take the partial derivatives of $g(y,\nu)$. Only the first three derivatives in terms of $y$, the first derivative with respect to $\nu$, and the first mixed derivative are required in this case. By using the Taylor expansion of this system, we obtain the following derivatives:
\begin{equation}\label{eq:g}
\begin{array}{rcl}
g_y &=& \langle e, d\Phi(0,\nu_{k}^{*})(e + w_y)\rangle \\
g_{y^2} &=& \langle e, d\Phi(0,\nu_{k}^{*})(w_{y^2})+ d^2\Phi(0,\nu_{k}^{*})(e+w_y,e+w_y)\rangle\\
g_{y^3} &=& \langle e, d\Phi(0,\nu_{k}^{*})(w_{y^3})+ 3d^2\Phi(0,\nu_{k}^{*})(e+w_y,w_{y^2})+d^3\Phi(0,\nu_{k}^{*})(e+w_y,e+w_y,e+w_y)\rangle\\
g_{\nu} &=& \langle e, d\Phi(0,\nu_{k}^{*})( w_{\nu})+\Phi_{\nu}(0,\nu_{k}^{*})\rangle \\
g_{y\nu} &=& \langle e, d\Phi_{\nu}(0,\nu_{k}^{*})(e+w_{y})+d\Phi(0,\nu_{k}^{*})(w_{y\nu}) +d^2\Phi(0,\nu_{k}^{*})(e+w_y, w_{\nu})\rangle
\end{array}
\end{equation}
We obtain the explicit values of the first few derivatives of $w$ with respect to $y$ and $\nu$ from the equation $G(y,w(y,\nu),\nu)\equiv 0$, those are:
\begin{equation}\label{eq:solve}
\begin{array}{rcl}
0&=&Pd\Phi(0,\nu_{k}^{*})(e + w_y) \\
0&=&Pd\Phi(0,\nu_{k}^{*})(w_{y^2})+ Pd^2\Phi(0,\nu_{k}^{*})(e+w_y,e+w_y)\\
0&=&Pd\Phi(0,\nu_{k}^{*})(w_{y^3})+ P3d^2\Phi(0,\nu_{k}^{*})(e+w_y,w_{y^2})+Pd^3\Phi(0,\nu_{k}^{*})(e+w_y,e+w_y,e+w_y)\\
0&=&Pd\Phi(0,\nu_{k}^{*})( w_{\nu})+P\Phi_{\nu}(0,\nu_{k}^{*})\\
0&=&Pd\Phi_{\nu}(0,\nu_{k}^{*})(e+w_{y})+Pd\Phi(0,\nu_{k}^{*})(w_{y\nu}) +Pd^2\Phi(0,\nu_{k}^{*})(e+w_y, w_{\nu})
\end{array}
\end{equation}
For both sets of equations~\eqref{eq:g} and~\eqref{eq:solve}, the derivatives of $w$ are evaluated at $(0,\nu_{k}^{*})$ and we now obtain explicit expressions for them. Recall that $L_0 = d\Phi(0,\nu_{k}^{*})$ and $PL_0=L_0$. We begin with
\[
0=Pd\Phi(0,\nu_{k}^{*})(e + w_y) = L_0 (e) + L_0(w_y) = L_0 w_y.
\]
But $w_y(0,\nu_{k}^{*}) \in (\ker L_0)^\bot$ because $w(y,\nu) \in (\ker L_0)^\bot$. By the invertibility of $L_0$ on $(\ker L_0)^{\bot}$ we have $w_y(0,\nu_{k}^{*})=0$ and so $g_y(0,\nu_{k}^{*})=0$ from the first equation in~\eqref{eq:g}.

From a Taylor expansion of equation (\ref{KS-op}), one can show that $d^2\Phi(0,\nu_{k}^{*})(\zeta_1,\zeta_2)=\frac{d}{dx}(\zeta_1\zeta_2)$
where $\zeta_1,\zeta_2$ are functions in the tangent space of $X$ at $(0,\nu_{k}^{*})$. This is applied to $Pd\Phi(0,\nu_{k}^{*})(w_{y^2})+ Pd^2\Phi(0,\nu_{k}^{*})(e+w_y,e+w_y)=0$. Because $w_y=0$, this implies
\[
PL_0 w_{y^2}(0,\nu_{k}^{*})+P[(e^2)']=0,
\]
where $(e^2)' = \left(\cos^2\left(\beta\xi\right)\right)'=-\beta \sin\left(2\beta\xi\right)$
from which we get
\[
w_{y^2}(0,\nu_{k}^{*})=\beta L_0^{-1} \sin\left(2\beta\xi\right)
=\dfrac{\beta}{(2\beta)^2(\nu_{k}^{*}(2\beta)^2-1)}\sin\left(2\beta\xi\right) = \dfrac{1}{6\pi (2n-1)}\sin(2\beta\xi).
\]
Note that $\sin$ is an eigenfunction of $L_0$ and so the $L_0^{-1}\sin(2\beta\xi)$ computation is obtained by inverting the corresponding eigenvalue.
Now, $g_{y^2} = \langle e, d\Phi(0,\nu_{k}^{*})(w_{y^2})+ d^2\Phi(0,\nu_{k}^{*})(e+w_y,e+w_y)\rangle$ and by substitution of both $L_0$ and $(e^2)'$ we obtain
\[
g_{y^2}(0,\nu_{k}^{*})=\langle e,L_0w_{y^2}(0,\nu_{k}^{*}) + (e^2)' \rangle = 0
\]
because $\langle e, (e^2)' \rangle=0$ by oddness of the integrand and $\langle e, L_0 w_{y^2} \rangle$ vanishes by orthogonality of
$\cos$ and $\sin$.
Notice that in the expression for $g_{y_3}$ from~\eqref{eq:g}, we have $\langle e,L_0 (w_{y^3})\rangle = 0$ by orthogonality of $\ker L_0$ with
$\mbox{range}\,L_0$. By Taylor expansion of $\Phi$ we have $d^3 \Phi=0$ so we have the simplification to
\[
g_{y^3}(0,\nu_{k}^{*})=\langle e, 3(ew_{y^2}(0,\nu_{k}^{*}))' \rangle =  \left\langle e, \frac{1}{4}(\cos(3\beta\xi)+\cos(\beta\xi)\cos(2\beta\xi)) \right\rangle = \frac{1}{8}.
\]
Since $G(y,w(y,\nu),\nu)\equiv 0$, then $P\Phi_{\nu}(0,\nu_{k}^{*})=0$ which means the last equality of~\eqref{eq:solve} yields $w_{\nu}(0,\nu_{k}^{*})=0$, thus $g_{\nu}(0,\nu_{k}^{*}) =0$. Finally, $d\Phi_{\nu}\xi  =\xi^{''''}$, implies $Pd\Phi_{\nu}(0,\nu_{k}^{*})(e)=\beta^4 P(e)=0$ and $d^{2}\Phi(0,\nu_{k}^{*})(e,0)=0$ by the formula for the second derivative. This leads to $w_{y\nu}(0,\nu_{k}^{*})=0$ and we obtain
\[
g_{y\nu}(0,\nu_{k}^{*})=\langle e, e'' \rangle = -\beta^2.
\]
Thus, we have $g(0,\nu_k^{*})=g_{y}(0,\nu_k^{*})=g_{yy}(0,\nu_k^{*})=g_{\nu}(0,\nu_k^{*})=0$ and $g_{yyy}(0,\nu_k^{*})g_{y\nu}(0,\nu_k^{*})<0$.
By standard results from singularity theory this indicates a supercritical pitchfork bifurcation \cite{GS85}. For the bifurcation point at
$\nu_{1/2}^{*}$, this automatically guarantees that the bifurcating branches of steady-states are asymptotically stable.
\qed

\section{Numerical continuation\label{numerical}}

\begin{figure}[h]
\begin{center}
\includegraphics[width=0.8\textwidth]{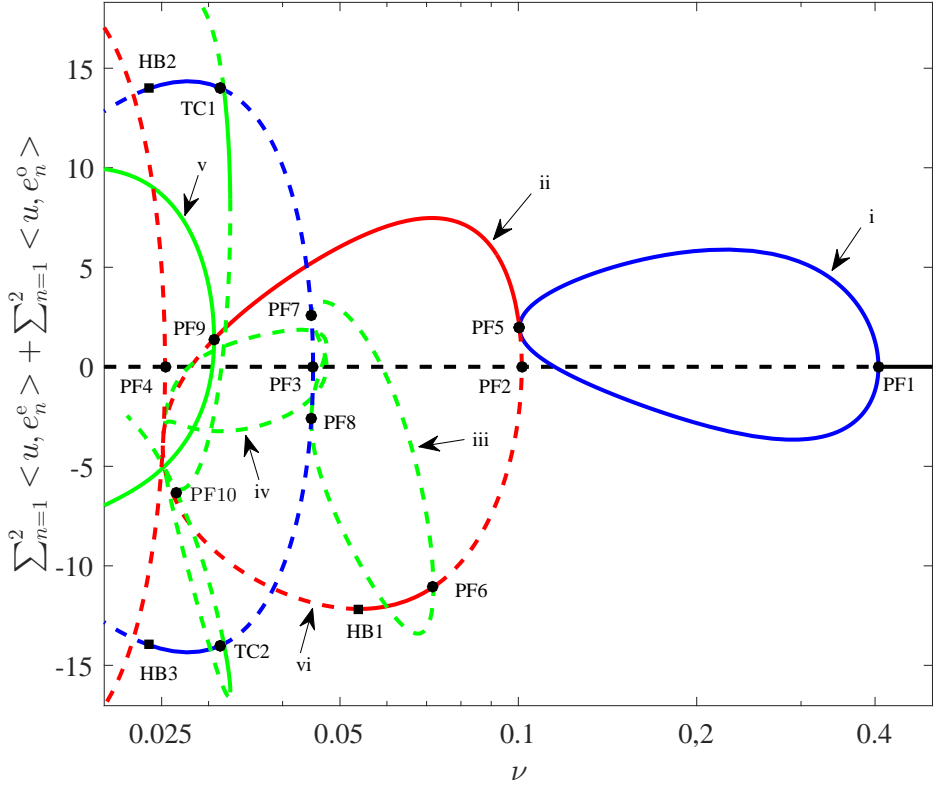}
\end{center}
\caption{Partial bifurcation diagram of the KS IBVP for homogeneous Dirichlet boundary conditions. Shown are the primary branches emanating
at $\nu=1/(\pi k)^2$ for $k=1/2, 3/2$ in blue and for $k=1,2$ in red, and secondary branches
bifurcating off the latter in green. Solid lines denote stable branches and dashed lines denote unstable branches.
Pitchfork bifurcations are denoted by solid circles and Hopf bifurcations by solid squares. On the vertical axis the
sum of projections onto the kernel functions $e^{\rm e}_n$ and $e^{\rm o}_n$ ($n=1,2$) is shown. The Roman numerals {\bf i}--{\bf vi} correspond
to the solutions visualized in Figure~\ref{portraits}.}
\label{KS_bif1}
\end{figure}

Having determined that branches of equilibria with (hidden) symmetry emanate from the zero solution in supercritical pitchfork bifurcations,
we turn our attention to the global behaviour of these solutions.
Diagram~\ref{KS_bif1} shows the primary branches with wave number $k$ equal to $1/2$ and $3/2$ in blue and
to $1$ and $2$ in red. All bifurcation points at which primary branches lose -- or gain -- stability are included and the corresponding
secondary branches are shown in green. Stable solutions are shown with solid lines and unstable solutions with dashed lines.
The upper and lower parts of the primary branches
corresponding to wave numbers $1/2$ and $3/2$ are related by the reflection symmetry $\kappa$ and thus display the same bifurcations.

The branch with $k=1/2$ is entirely stable and connected to the branch with $k=1$ at pitchfork bifurcation $\text{PF}_5$. The top branch of the $k=1$ family turns stable here, and spawns a stable branch of equilibria without reflection symmetry in $\text{PF}_9$. The bottom branch of this family also turns stable in a pitchfork bifurcation, $\text{PF}_6$, and subsequently produces a stable time-periodic solution in Hopf bifurcation $\text{HB}_1$. The branch with $k=3/2$ turns
stable in a transcritical bifurcation, as expected given the symmetry considerations explained in Section~\ref{sec:hidden}, and then exhibits the second Hopf bifurcation
around $\nu=0.025$. Around this parameter value, two symmetry-related equilibria on the above-mentioned stable secondary branch co-exist with two
symmetry-related stable equilibria on the $k=3/2$ branch, or two stable time-periodic solutions. The primary branch from $\nu_k^{*}$ with $k=2$ is unstable down 
to $\nu=0.02$, where we stopped the continuations.

\begin{figure}[p]
\begin{center}
\includegraphics[width=0.7\textwidth]{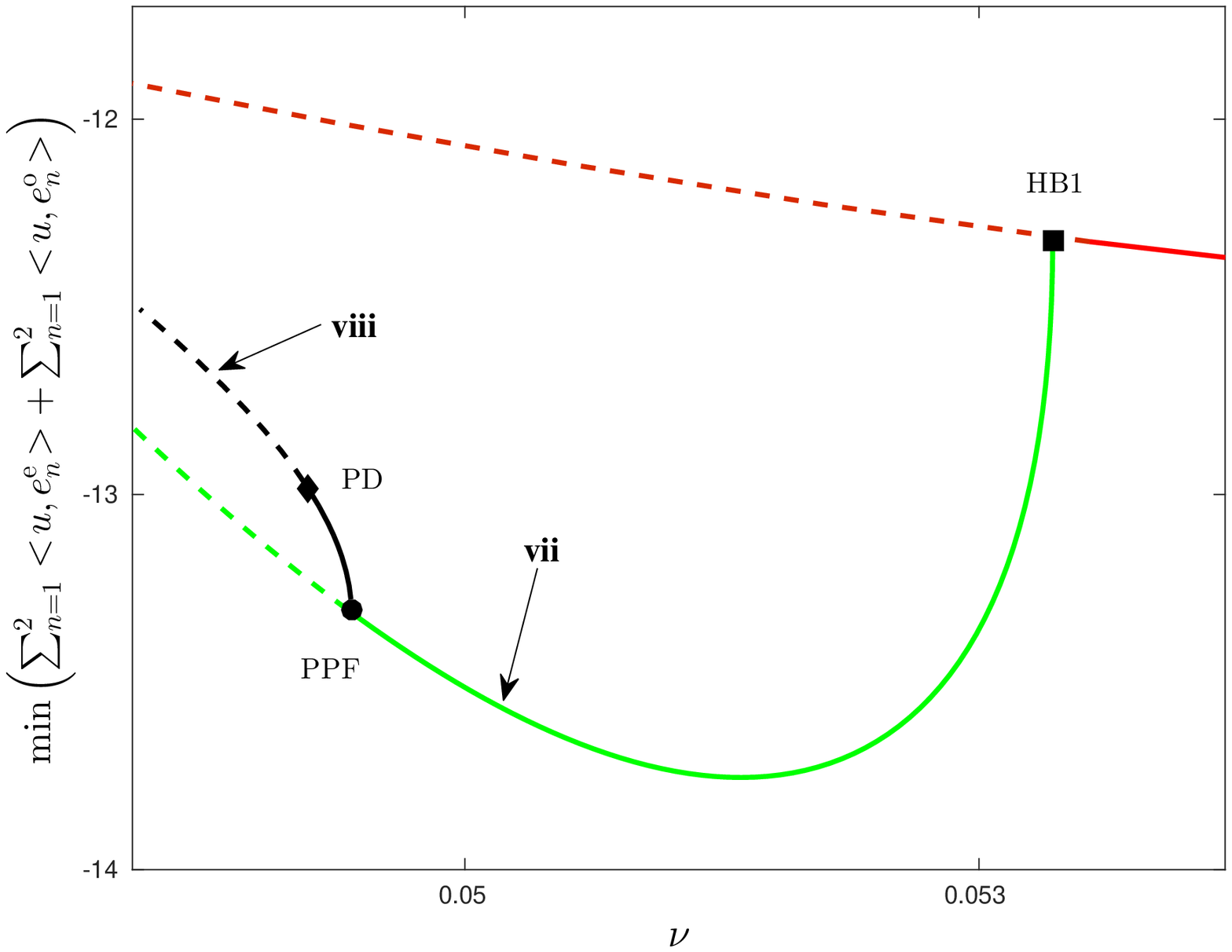}
\end{center}
\caption{Detail of diagram \ref{KS_bif1} around HB1. PPF denotes a pitchfork birucation of a periodic orbit and the solid diamond denotes
a period doubling bifurcation. The roman numerals {\bf vii}--{\bf viii} correspond
to the solutions visualized in Figure~\ref{portraits}.}
\label{KS_bif2}
\end{figure}

\begin{figure}[p]
\begin{center}
\includegraphics[width=0.7\textwidth]{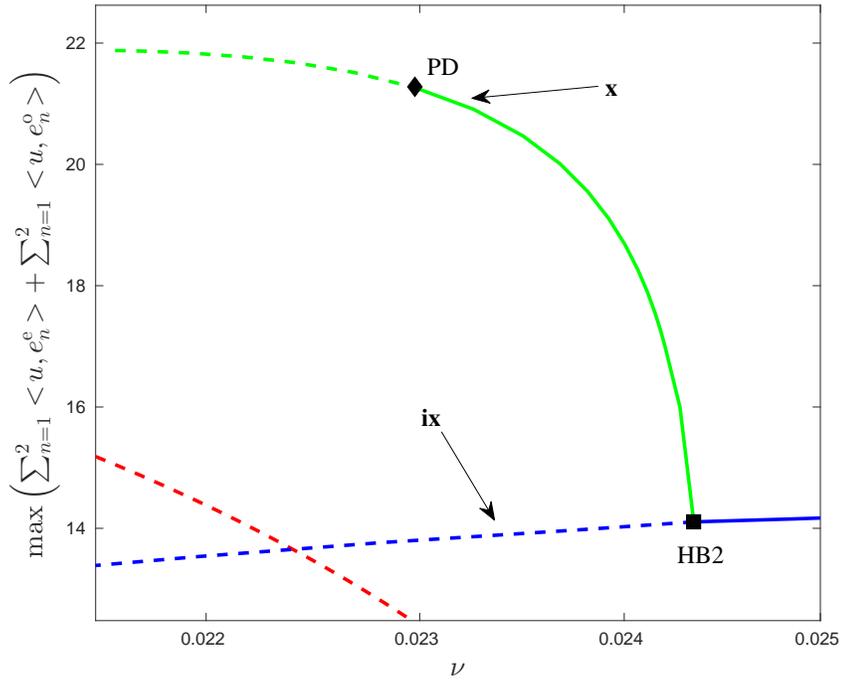}
\end{center}
\caption{Detail of diagram \ref{KS_bif1} around HB2. The solid diamond denotes
a period doubling bifurcation. The roman numerals {\bf ix}--{\bf x} correspond
to the solutions visualized in Figure~\ref{portraits}.}
\label{KS_bif3}
\end{figure}

Details of the bifurcation diagram around the Hopf bifurcations are shown in Figures~\ref{KS_bif2} and~\ref{KS_bif3}.
The periodic solutions originating from $\text{HB}_1$ inherit a discrete symmetry from the $k=1$ branch, namely a shift over half the period
followed by reflection $\kappa$. This symmery is broken in pitchfork bifurcation $\text{PPF}$. The branch of periodic solutions from
$\text{HB}_2$ has no symmetry. Both branches of periodic solutions undergo period doubling bifurcations that could signal the onset of more
complicated spatio-temporal behaviour, beyond the scope of the current paper.

\begin{figure}[h]
\begin{center}
\includegraphics[height=0.5\textheight]{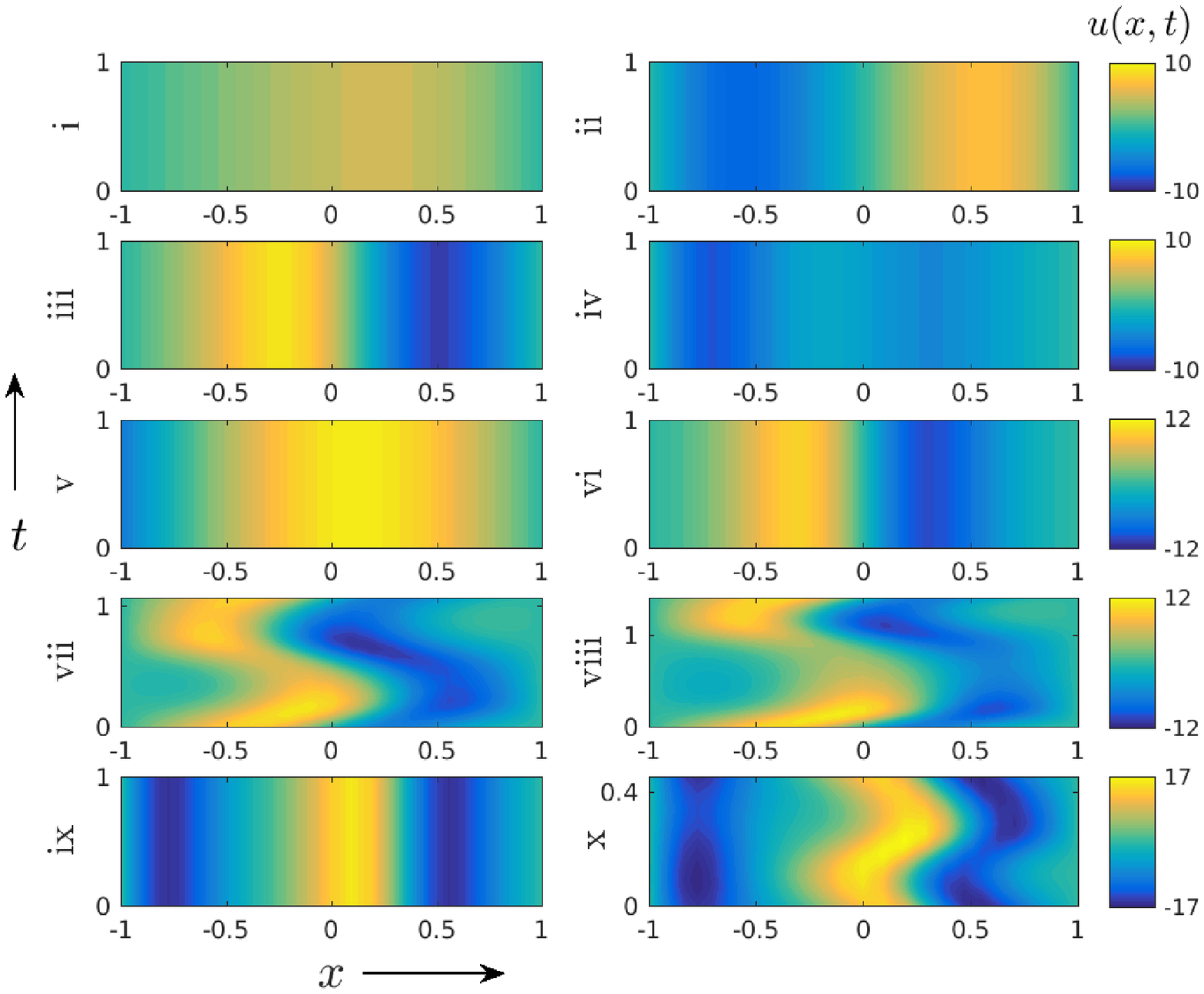}
\end{center}
\caption{Some solutions corresponding to the labels in Figures~\ref{KS_bif1}--\ref{KS_bif3}. The horizontal axes denote space and the
vertical axes denote time. The solutions $u(x,t)$ are shown using the colour scale on the right hand side of each row.}
\label{portraits}
\end{figure}

Along the various branches of equilibria and periodic orbits we have visualized a number of solutions. They are shown in Figure~\ref{portraits}. 
The spatial structure of the solutions at these relatively large values of the viscosity is dominated by wave numbers up to two, and the amplitude 
increases steadily with decreasing viscosity.

\section{Conclusion}

Using Lyapunov-Schmidt reduction and numerical continuation, and exploiting the hidden symmetry, we have given a fairly complete overview
of the bifurcation structure of the KS equation with Dirichlet BC near the onset of nontrivial behaviour. In confluence with earlier
results by \citet{LiChen}, this overview includes a complete description of all bifurcations from the
zero solution. These are all supercritical pitchfork bifurcations, but come in two families: one breaks the reflection symmetry of the
system under consideration, while the other breaks a reflection symmetry of an extended problem with periodic boundary conditions and is
thus due to hidden symmetry.

Using numerical continuation, we computed the first four primary branches, the bifurcation points at which these
lose stability and the resulting secondary branches. Steady-state bifurcations from the primary branches can be of the pitchfork and transcritical
type, and both are explained in terms of hidden symmetry. The secondary branches include time-periodic solutions that exhibit
pitchfork as well as period doubling bifurcations. The resulting partial bifurcation diagram, shown in Figure~\ref{KS_bif1},
is as complicated as one might expect from such a highly nonlinear system, and includes a wide range of parameters in which
several stable equilibria and periodic solutions co-exist.

In addition to elucidating the transition from trivial to spatio-temporally structured behaviour in the KS problem with Dirichlet boundary conditions, we
have, with this study, provided a neat example of the consequences of hidden symmetry. This phenomenon has mostly been
demonstrated through evolution equations of the reaction-diffusion type, for instance by \citet{Craw} and~\citet{GolSte}. We hope that the KS example
will prove a useful addition to literature on hidden symmetry and aid in future studies of parabolic initial-boundary value problems.

\nonumsection{Acknowledgments} \noindent EF was supported by an Ontario Graduate Scholarship, PLB and LvV are
supported by the National Science and Engineering Research Council of Canada in the form of a Discovery Grant.

\appendix{Numerical simulation and continuation}

The continuation and stability analysis presented in Section~\ref{numerical} is based on the numerical simulation of IBVP (\ref{KS}--\ref{KS-DBC}).
Details of the discretization scheme can be found in \citet{LvV}, here we present only its extention to the linearized
equations and the application of Newton-Krylov and Arnoldi iteration.

The simulation algorithm is based on that of \citet{Rothe}: first, time
is discretized and then a linear BVP is solved for each time step with the aid of Green's function. Since we only need to integrate
the equations over short time intervals for the purpose of continuation of equilibria and periodic orbits, the semi-implicit Euler
method is adequate for time discretization. The BVP to be solved in each time step is then
\begin{align}
\mathcal{L}u^{(k+1)}\equiv(1+h\partial^2_x+h\nu\partial^4_x)u^{(k+1)} =u^{(k)}- \frac{h}{2} (u^{(k)})^2_x \\
u^{(k+1)}(-1)=u^{(k+1)}_{xx}(-1)=u^{(k+1)}(1)=u^{(k+1)}_{xx}(1)&=0.
\end{align}
and perturbing $u$ by $w$ and $\nu$ by $\omega$, we find the linearized equation
\begin{align}
\mathcal{L}w^{(k+1)}=w^{(k)}- h (u^{(k)}w^{(k)})_x -h\omega u^{(k+1)}_{xxxx}\\
w^{(k+1)}(-1)=w^{(k+1)}_{xx}(-1)=w^{(k+1)}(1)=w^{(k+1)}_{xx}(1)&=0.
\end{align}
where $h$ is the time step size. The solution is given in terms of Green's function for the differential operator $\mathcal{L}$ with the given boundary conditions as
\begin{align}
u^{(k+1)}&=G\ast u^{(k)} +\frac{h}{2} DG \ast (u^{(k)})^2 \\
w^{(k+1)}&=G\ast \left(w^{(k)}-h\omega u^{(k+1)}_{xxxx}\right) +h DG \ast (u^{(k)}w^{(k)})
\end{align}
where the asterisk denotes the spatial convolution.
Finally, the fourth derivative of the solution is time stepped according to
\begin{equation}
u_{xxxx}^{(k+1)}=D^4 G\ast u^{(k)} +\frac{h}{2} D^5 G \ast (u^{(k)})^2
\end{equation}
Green's function is known analytically, so the only numerical approximation is a quadrature rule for approximating
the convolutions. Following the approach of \citet{LvV}, we represent the solution on a global, closed Chebyshev grid
and employ an exponentially accurate combination of barycentric interpolation and Clenshaw-Curtis quadrature.
For $\nu>0.03$, we used a grid order of $N=32$ and below that we set $N=64$, while the time step was
fixed to $h=10^{-3}$. Several solutions
along each of the curves were recomputed at double their resolution and smaller time step for verification.

Following the approach explained in detail by \citet{sanch}, the time-steppers for the IBVP and its linearization are all that is needed
for the continuation and stability analysis of equilibria and periodic orbits. Hence, we need only state briefly the continuation equations, tolerances
and convergence properties of the algorithm.

The continuation equations are best formulated in terms of the flow of the vector of grid point values $\bm{u}\in\mathbb{R}^{N+1}$,
denoted by $\phi(\bm{u},t,\nu)$.
The equations to solve for equilibria and periodic orbits are
\begin{equation}\label{PAC1}
\bm{\phi}(\bm{u},P,\nu)-\bm{u}=0
\end{equation}
and a phase condition $\psi(\bm{u},P)=0$, namely
\begin{alignat}{3}\label{PAC2}
\psi\equiv P-c&=0\ \text{(for equilibria), or}&\quad \psi\equiv u_{n/2}&=0\ \text{(for periodic orbits)}
\end{alignat}
where $c=1$ is an integration time long enough to pre-condition the linear system for the
Newton-Raphson update step. We compute branches of solutions to BVP (\ref{PAC1}--\ref{PAC2}) using pseudo-arclength
continuation. The residual thresold is fixed to $10^{-8}$ in the norm
$$
\|(\bm{u},P,\nu)\|=\sqrt{\|\bm{u}\|_2^2+\psi^2}
$$
where $\|.\|_2$ denotes the approximate $L_2$--norm of a function as obtained by  Clenshaw-Curtis quadrature on the
grid. We obtain quadratic convergence of the Newton-Raphson iterates with about $10$ Krylov subspace iterations
and a tolerance of $10^{-6}$ for the linear solving.

\bibliographystyle{ws-ijbc}

\end{document}